\documentclass{amsart}

\usepackage{amssymb}

\newtheorem{thm}{Theorem}[section]

\newtheorem{lem}[thm]{Lemma}
\newtheorem{cor}[thm]{Corollary}
\theoremstyle{definition}

\theoremstyle{remark}
\newtheorem{remark}[thm]{Remark}
\numberwithin{equation}{section}
\begin{document}

\title[Lowest Eigenvalues]{Lowest eigenvalues and formally self-adjoint fourth order elliptic differential operators}

\author{David Raske}

\address{David Raske, 3635 Washtenaw Ave, Ann Arbor, Michigan, 48104, United States}

\email{nonlinear.problem.solver@gmail.com}


\begin{abstract}

Let $(M,g)$ be a closed, smooth, Riemannian manifold of dimension $m \geq 1$. Let $\eta$ be a smooth $(0,1)$-tensor field on $M$. The divergence of $\eta$ is defined as $\text{div}_g(\eta):=g^{ij}(\nabla \eta)_{ij}$. Now let $\Delta_g$ be a differential operator on $M$ that is given on functions by $\Delta_g u  = \text{div}_g \nabla u$. We will call $\Delta_g$ the Laplace-Beltrami operator. With this definition in place, it is not difficult to produce an example of a formally self-adjoint elliptic differential operator on $M$ that has a  sign-changing eigenfunction that is associated with the operator's lowest eigenvalue. Indeed, let $\lambda_2$ be the second lowest eigenvalue of $-\Delta_g$, and let $L_g$ be a differential operator on  $M$ that is given on functions by $L_g u = \Delta^2_g u + \lambda_2 \Delta u$. Then $L_g$ will possess a sign-changing eigenfunction that is associated with $L_g$'s lowest eigenvalue.. The question that remains is given a smooth, closed manifold $M$ of dimension $m \geq 1$, how rare are formally self-adjoint elliptic differential operators on $M$ that have sign-changing eigenfunctions that are associated with the operators' lowest eigenvalues. In this paper, we will see that if $A = T -\lambda g$, where $T$ is a smooth, symmetric, negative semi-definite $(0,2)$-tensor field on $M$, then $P_g$, the differential operator on $M$ given on functions by $P_gu=,\Delta_g^2 - \text{div}_g(A(\nabla u)^\sharp)$, will have the property that it possesses a sign-changing eigenfunction that is associated with the lowest eigenvalue of the operator. This suggests that on any smooth, closed manifold of dimension $m \geq 1$ there exists a lot of formally self-adjoint fourth order elliptic differential operators on the manifold that possess sign-changing eigenfunctions that are associated with the lowest eigenvalues of the operators.  

\end{abstract} 

\maketitle

\section{Introduction}

First, we must define some objects. Let $(M,g)$ be a closed, smooth Riemannian manifold of dimension $m \geq 1$. Let $k$ be an even positive integer. Suppose $P$ is a $k$th order,  formally self-adjoint elliptic differential operator on $M$. We will call $P$ a kth order Type 1 formally self-adjoint elliptic differential operator on $M$ if and only if (a) there exists a real number $\lambda^*$ such that there exists a smooth, signed solution $u$ of the partial differential equation $Pu = \lambda^*u$; (b) if there exists another smooth, not identically zero solution $z$ of $Pz=\lambda^* z$, then there exists a real number $c$ not equal to zero such that $z = cu$; and (c) there does not exist a smooth, not identically zero solution of the partial differential equation $Pw=\lambda w$ if $\lambda < \lambda^*$.  If a kth order formally self-adjoint elliptic differential operator on $M$ is not a kth order Type 1 formally self-adjoint elliptic differential operator on $M,$ we will call it a kth order Type 2 formallly self-adjoint elliptic differential operator on $M$. For the sake of brevity, we will refer to kth order Type 1 formally self-adjoint elliptic differential operators on $M$ as kth order Type 1 FSA elliptic differential operators on $M$, and we will refer to kth order Type 2 formally self adjoint elliptic differential operators on $M$ as kth order Type 2 FSA elliptic differential operators on $M$.

The eigenproblem for second order elliptic differential operators is very well understood. Let $(M,g)$ be a closed, smooth Riemannian manifold of dimension $m \geq 1$. Then the Krein-Rutman Theorem  (see \cite{Sweers} for a discussion of this theorem) guarantees that if we let $h \in C^\infty(M)$, and if we let $C_g$ be a differential operator on $M$ that is given on functions by $C_g u= -\Delta_g u+ h u$, then $C_g$ is  a second order Type 1 FSA elliptic differential operator on $M$. \footnote{$\Delta_g$ is defined in the abstract.} The situation is very different for higher-order elliptic differential operators.  Indeed, let $\lambda_2$ be the second lowest eigenvalue of $-\Delta_g$. One of the simplest fourth order elliptic differential operators on $M$ that one could consider, $L_g$, which is  given on functions by $L_g u= \Delta_g^2 u + \lambda_2 \Delta_g u$, has the property that there exists a sign-changing $C^\infty(M)$ function $w$ such that $L_g w =\lambda w$, where $\lambda$ is the lowest eigenvalue of the operator $L_g$. (We will prove this in Section Two.) Given this, it is natural to ask how common are fourth order Type 2 FSA elliptic differential operators on M? Let $a$ and $b$ be positive real numbers. The example that one typically finds are differential operators on $M$ that are of the form $\Delta^2_g + a\Delta_g + b$. (See \cite{HR} and \cite{Robert} for more on these operators.) In this paper we will extend our focus to differential operators on $M$ of the form $\Delta_g^2 - \text{div}_g(A(\nabla \cdot)^\sharp) + b$, where $A$ is a smooth, symmetric $(0,2)$ tensor field on $M$ and where $b$ is a real number. Notice that this kind of operator includes operators of the form $\Delta_g^2 + a \Delta_g + b$. Differential operators of the form $\Delta^2_g - \text{div}_g(A(\nabla \cdot)^\sharp) + b$ also occur naturally in the study of prescription of $Q$-curvature on closed, smooth Riemannian manifolds of dimension greater than four. (See \cite{Robert} for more on the fascinating subject.) After some investigation we arrive at the following:

\begin{thm} Let $(M,g)$ be a smooth, closed Riemannian manifold of dimension $m \geq 1$. Let $\lambda_2$ be the second lowest eigenvalue of $-\Delta_g$.  Let $A$ be a smooth, symmetric $(0,2)$-tensor field on $M$. Let $P_g$ be a fourth order differential operator on $M$ that is given on functions by $P_g u = \Delta^2_g u- \text{div}_g(A(\nabla u)^\sharp)$. Then, $P_g$ has a sign-changing  eigenfunction that is associated with the lowest eigenvalue of $P_g$, if $A$ has the property that $\max_{x \in M} \gamma_m(x) \leq -\lambda_2$ where $\gamma_m(x)$ is 
$$
\max_{w \in T_x(M) \setminus \{0\}}  \frac{A(x)(w,w)}{g(x)(w,w)}.
$$
($A(x)$ is  $A$ at $x \in M$, and $g(x)$ is $g$ at $x \in M$.) This, in turn, allows us to conclude that if $A$ has the property that $\max_{x \in M} \gamma_m(x) \leq -\lambda_2$, then $P_g$ is a fourth order Type 2 FSA elliptic differential operator on $M$.
\end{thm}

Now recalling some basic algebra, we can conclude the following

\begin{cor} Let $(M,g)$ be a smooth, closed Riemannian manifold of dimension $m \geq 1$. Let $\lambda_2$ be the second lowest eigenvalue of $-\Delta_g$.  Let $A$ be a smooth, symmetric $(0,2)$-tensor field on $M$. Let $P_g$ be a fourth order differential operator on $M$ that is given on functions by $P_g u = \Delta^2_g u - \text{div}_g(A(\nabla u)^\sharp)$.
Then, $P_g$ has a sign-changing  eigenfunction that is associated with the lowest eigenvalue of $P_g$, if there exists a smooth, symmetric, negative semi-definite $(0,2)$ tensor field $T$, such that $A= T-\lambda_2g$. This, in turn, allows us to conclude that if there exists a smooth, symmetric negative semi-definite $(0,2)$-tensor field  $T$ such that $A=T-\lambda_2 g$, then $P_g$ is a fourth order Type 2 FSA elliptic differential operator on $M$.
\end{cor}

\begin{remark} Let $(M,g)$ be a closed, smooth Riemannian manifold of dimension $m \geq 1$. Let $k$ be a positive even integer. Let $c$ be a real number, and let $A$ be a smooth, symmetric $(0,2)$-tensor field on $M$. Let $P_g$ be a differential operator on $M$ that is given on functions by $P_g u = \Delta^2_g u - \text{div}_g(A(\nabla u)^\sharp)$. Let $F_p$ be a differential operator given on functions by $F_g u = \Delta^2_g u- \text{div}_g(A(\nabla u)^\sharp) +cu$. If $P_g$ is a fourth order Type 2 FSA elliptic differential operator on $M$, then $F_g$ is a fourth order Type 2 FSA elliptic differential operator on $M$. 
\end{remark}

Let $(M,g)$ be a smooth, closed Riemannian manifold of dimension $m \geq 1$.The proof of Theorem 1.1 relies heavily on the variational characterization of the lowest eigenvalue of a formally self-adjoint elliptic differential operator on $M$ and the fact that constants will always be in the kernel of the operator $P_g$, where $A$ is a smooth, symmetric $(0,2)$-tensor field on $M$ and $P_g$ is the differential operator on $M$ given on functions by $P_g u = \Delta_g^2 u-\text{div}_g(A(\nabla u)^\sharp)$. Noting these two facts, we see that the operator $P_g$ will be a fourth order Type 2 FSA elliptic differential operator on $M$ if there exists a non-constant smooth function $u$ such that $\int_M u P_g u \, \text{dvol}_g \leq 0$. We then see that if $A$ is "negative enough", then $\phi$, an eigenfunction associated with the second lowest eigenvalue of the elliptic differential operator on $M$, $-\Delta_g$, will have the property that $\int_M \phi P_g \phi \, \text{dvol}_g \leq 0.$ The proof of Theorem 1.1 can be found in the next section. The lemmas used to prove Theorem 1.1 can be found there as well.

\section{Results and Proofs}
First we need to define some Hilbert spaces. The Hilbert space $L^2(M)$ is the set of locally integrable functions
$u$ on $M$ for which the norm
$$
||u||_2 = (\int_M u^2 \, \text{dvol}_g)^{1/2}
$$
is finite. The $L^2(M)$ inner-product $<\cdot,\cdot>_2$, which is a real-valued function defined on $L^2(M) \times L^2(M)$,  is given by $<u,v>_2 = \int_M u v \, \text{dvol}_g$, for all $u,v \in L^2(M)$.

If $k$ is a positive integer, the Hilbert space $H^{k}(M)$ is the set of $u \in L^2(M)$ such that $Pu =f \in L^2(M)$ (in the weak sense) whenever $P$ is a differential operator of order $\leq k$. We define the Hilbert space norm $|| \,||_{H^k(M)}$ on $H^k(M)$ by
$$
||u||_{H^k(M)} = (\Sigma_{i=0}^k \int_M |\nabla^i u|^2 \text{dvol}_g)^{1/2}.
$$
$H^k(M)$ comes equipped with an inner-product, $<\cdot,\cdot>_{H^k(M)}$, which is a real valued function defined on $H^k(M) \times H^k(M)$, such that 
$$
<u,u>_{H^k(M)} = ||u||_{H^k(M)}^2
$$
for all $u \in H^k(M)$. We won't be using this inner-product in the sequel, so we won't write it down here. We will use the fact that the norm $||\,||_{H^2(M)}$ is equivalent to the norm $||\,||_{H^2_*(M)}$, where $||u||_{H^2_*(M)} = (\int_M ((\Delta_g u)^2 + u^2) \, \text{dvol}_g)^{1/2}$, for all $u \in H^2(M)$..

The first lemma looks simple but will prove to be very helpful later on. It is as follows.

\begin{lem}
Let $(M,g)$ be a closed, smooth Riemannian manifold of dimension $m \geq 1$. Let $u \in H^2(M)$. Then we have
$$
\int_M g((\nabla u)^\sharp,(\nabla u)^\sharp) \, \text{dvol}_g \leq (\int_M (\Delta_g u)^2 \, \text{dvol}_g)^{1/2} (\int_M u^2 \, \text{dvol}_g)^{1/2}.
$$ 
\end{lem} 

\begin{proof} First note that we have 
$$
\int_M g((\nabla u)^\sharp,(\nabla u)^\sharp) \, \text{dvol}_g = \int_M (-\Delta_g u)(u) \, \text{dvol}_g
$$
 for all $u \in C^\infty(M)$. Applying H{\"o}lder's inequality to the right hand side of this equation, we see that
$$
\int_M (-\Delta_g u)(u) \, \text{dvol}_g \leq (\int_M (\Delta_g u)^2 \, \text{dvol}_g)^{1/2} (\int_M u^2 \, \text{dvol}_g)^{1/2}
$$
for all $u \in C^\infty(M)$. Recalling that $C^\infty(M)$ is dense in $H^2(M)$, we have the lemma.
\end{proof}

We can now turn our attention to proving the most important lemma of the paper.

\begin{lem}
Let $(M,g)$ be a smooth, closed Riemannian manifold of dimension $m \geq 1$. Let $A$ be a smooth, symmetric $(0,2)$-tensor field on $M$. Let $P_g$ be the differential operator on $M$ that is given on functions by $P_g u = \Delta^2_g  u - \text{div}_g(A(\nabla u)^\sharp)$. Let $F: H^2(M) \setminus \{0\}\rightarrow \mathbb{R}$ be defined as follows: 
\begin{equation}
F(u) = \frac{\int_M ((\Delta_g u)^2 + A((\nabla u)^\sharp ,(\nabla u)^\sharp) \, \text{dvol}_g}{\int_M u^2 \, \text{dvol}_g},
\end{equation}
Then 
$$
\inf_{H^2(M) \setminus \{0\}} F(u)
$$ 
exists (let's call it $\lambda$), and $F$ is minimized over $H^2(M) \setminus \{0\}$ by $C^\infty(M)$ functions that satisfy the equation $P_g w = \lambda w$. We also have that $\lambda$ is the lowest eigenvalue of $P_g$.

\end{lem}

\begin{proof}
First note that if  $\inf_{u \in H^2(M) \setminus \{0\}} F(u) = -\infty$, we have a sequence of $H^2(M)$ functions $\{u_n\}_{n=1}^\infty$ with 
$$
\int u_n^2 \, \text{dvol}_g = 1
$$ 
such that $\int A((\nabla u_n)^\sharp, (\nabla u_n)^\sharp)  \, \text{dvol}_g \rightarrow - \infty$ as $n \rightarrow \infty$. This, in turn, implies 
$$
\int g((\nabla u_n)^\sharp,(\nabla u_n)^\sharp) \,  \text{dvol}_g \rightarrow \infty 
$$
as $n \rightarrow \infty$. This can't be the case, though, because of Lemma 2.1. So we have it that the infimum of $F$ over $H^2(M) \setminus \{0\}$ is finite. Let's call it $\lambda$. 

We will now see that a minimizer of $F$ over $H^2(M) \setminus \{0\}$ exists. Let $\{u_n\}_{n=1}^\infty$ be a minimizing sequence, where $u_n \in H^2(M)$ for all positive integers $n$. Note that we can assume without loss of generality that $\int_M u_n^2 \, \text{dvol}_g = 1$ for all positive integers $n$. Now suppose $\int_M ((\Delta_g u_n)^2 + u_n^2) ,\ \text{dvol}_g \rightarrow \infty$ as $n \rightarrow \infty$.  Suppose as well that there exists a real number $B$ such that $F(u_n) \leq B$ for all $n \in \mathbb{N}$. It would follow then that $\int_M A((\nabla u_n)^\sharp,(\nabla u_n)^\sharp) \, \text{dvol}_g \rightarrow -\infty$ as $n \rightarrow \infty$. This, in turn, would imply that $\int_M g((\nabla u_n)^\sharp,(\nabla u_n)^\sharp) \, \text{dvol}_g \rightarrow \infty$ as $n \rightarrow \infty$. This contradicts Lemma 2.1, though, so we have it that $F(u_n) \rightarrow +\infty$ as $n \rightarrow \infty$ if  $\int_M ((\Delta_g u_n)^2 + u_n^2) ,\ \text{dvol}_g \rightarrow \infty$ as $n \rightarrow \infty$.  

We now have it that there exists a real number $D$ such that $||u_n||_{H^2(M)} \leq D$ for all positive integers $N$. This in turn implies that there exists a $H^2(M)$ function $z$ and a subsequence $\{u_{n_k}\}_{k=1}^\infty$ of $\{u_n\}_{n=1}^\infty$ that converges weakly with respect to the $H^2(M)$ norm to $z$ as $n \rightarrow \infty$. We wll now relabel the subsequence $\{u_{n_k}\}_{n=1}^\infty$ so that it becomes $\{ u_n\}_{n=1}^\infty$. This, in turn, allows us to write $u_n \rightarrow z$ with respect to the $H^1(M)$ norm, as $n \rightarrow \infty$. Now note that we can write 
\begin{equation}
\begin{split}
F(u_n) = &\int_M (((\Delta_g u_n)^2 + L((\nabla u_n)^\sharp, (\nabla u_n)^\sharp) \\ &+ \gamma u_n^2) + (- c g((\nabla u_n)^\sharp,(\nabla u_n)^\sharp) - \gamma u_n^2)) \, \text{dvol}_g.
\end{split}
\end{equation}
Here $L$ is a smooth, symmetric, positive definite $(0,2)$-tensor field on $(M,g)$, $\gamma$ is a positive real number, and $c$ is a positive real number. Now let $I: H^2(M) \rightarrow \mathbb{R}$, be defined as follows:
$$
 I(u) = \int_M  ((\Delta_g u)^2 + L((\nabla u)^\sharp,(\nabla u)^\sharp) + \gamma u^2) \, \text{dvol}_g
$$ 
for all $u \in H^2(M)$, Note that $I(\cdot)$ is the square of a norm for $H^2(M)$. Since norms are weakly lower semicontinuous we have it that the 
\begin{equation}
\liminf_{n \rightarrow \infty} I(u_n) \geq I(z).
\end{equation} 

Since we also have it that $u_n \rightarrow z$ with respect to the $H^1(M)$ topology, as $n \rightarrow \infty$, can write 
 $$
\int_M (-c g((\nabla u_n)^\sharp, (\nabla u_n)^\sharp)- \gamma u_n^2) \, \text{dvol}_g \rightarrow \int_M (-c g((\nabla z)^\sharp,(\nabla z)^\sharp) - \gamma z^2 ) \, \text{dvol}_g,
$$
as $n \rightarrow \infty$. It follows that $F(z) \leq \lambda$, and hence $F(z)= \lambda$, because $\lambda$ is the infimum of $F$ over $H^2(M) \setminus \{0\}$. (Note: we know that $z$ is not the zero function, because $\int_M z^2 \, \text{dvol} = 1$.)

Now let $w$ be a minimizer of $F$ over $H^2(M) \setminus \{0\}$. We then would have that $w$ satisfies the following equation:
\begin{equation}
\begin{split}
\int_M (\Delta_g w \Delta_g v +& A((\nabla w)^\sharp,(\nabla v)^\sharp)) \, \text{dvol}_g = \\ & ( \frac{\int_M ((\Delta_g w)^2 + A((\nabla w)^\sharp,(\nabla w)^\sharp)) \, \text{dvol}_g}{\int_M w^2 \, \text{dvol}_g}) \int_M w v \, \text{dvol}_g,
\end{split}
\end{equation}
for all $v \in H^2(M)$. We can now use elliptic regularity theory to see that $w$ is a smooth solution of the following equation: \footnote{Let $M$ be a smooth, closed Riemannian manifold of dimension $m \geq 1$, and let $A$ be a smooth, symmetric $(0,2)$-tensor field on $M$. On page 12 of \cite{Robert} it is shown that given $u,v \in C^\infty$ we have  $\int_M -\text{div}_g(A(\nabla u)^\sharp) v \, \text{dvol}_g = \int_M A((\nabla u)^\sharp,(\nabla v)^\sharp) \, \text{dvol}_g$.}
$$
P_g w=( \frac{\int_M ((\Delta_g w)^2 + A((\nabla w)^\sharp,(\nabla w)^\sharp)) \, \text{dvol}_g}{\int_M w^2 \, \text{dvol}_g})w. 
$$
It follows that $\lambda$ is an eigenvalue of $P_g$.

Now suppose that $P_g$ has an eigenvalue that is lower than $\lambda$. Let's call it $\eta$. This would imply that there exists a smooth, not identically zero function $y$ such that $P_g y = \eta y$, This, in turn would allow us to write
$$
 \frac{\int_M ((\Delta_g y)^2 + A((\nabla y)^\sharp,(\nabla y)^\sharp)\, \text{dvol}_g}{\int_M y^2 \, \text{dvol}_g} =  \frac{\int_M yP_gy \, \text{dvol}_g}{\int_M y^2 \, \text{dvol}_g} = \eta < \lambda.
$$
This contradicts the fact that $\lambda$ is the infimum of $F$ over $H^2(M) \setminus \{0\}$. The lemma follows.

\end{proof}

Finally, we prove a lemma that interacts with Lemma 2.2 nicely

\begin{lem} Let $(M,g)$ be a smooth, closed Riemannian manifold of dimension $m \geq 1$. Let $A$ be a smooth, symmetric $(2,0)$-tensor field on $M$. Let $P_g$ be the differential operator on $M$ that is given on functions by $P_g u = \Delta^2_g  u - \text{div}_g(A(\nabla u)^\sharp)$. If there exists a $C^\infty(M)$ function $u$ such that $G: H^2(M) \rightarrow \mathbb{R}$, where $G(z):=\int_M ((\Delta_g z)^2 + A((\nabla z)^\sharp, (\nabla z)^\sharp)$ has the property that $G(u) < 0$, then the lowest eigenvalue of the operator  $P_g$  is negative and the eigenfunctions associated with the lowest eigenvalue have mean values of zero.
\end{lem}

\begin{proof} First, recall that constants are in the kernel of $P_g$. This implies that $0$ is an eigenvalue of the differential operator $P_g$ and some of the eigenfunctions associated with that eigenvalue are constants. Note that we also have that $G(c)=0$ for all real numbers $c$. Now suppose that there exists a $H^2(M)$ function $u$ such that $G(u) < 0$. It follows from Lemma 2.2 that the functional $F$, which is defined in the statement of Lemma 2.2, is minimized by a smooth function $v$ and the minimum value, $\delta$, is negative. We can now use Lemma 2.2 to conclude that $\delta$ is the lowest eigenvalue of $P_g$. Now note that if $w: M \rightarrow \mathbb{R}$ is smooth, $\sigma$ is a real number, and $P_g w = \sigma w$, then we have that $0 =\int_M (P_g c) w \, \text{dvol}_g = \int_M c (P_g w) \, \text{dvol}_g = \sigma \int_M c w \, \text{dvol}_g$, where $c$ is a non-zero real number. It follows that if $\sigma$ is not $0$, the mean value of $w$ is zero. The lemma follows, because we have established that the lowest eigenvalue of $P_g$ is negative if there exists a function $u \in C^\infty(M)$ such that $G(u)<0$.
\end{proof} 

We can now prove Theorem 1.1, the main result of the paper.

\begin{thm} Let $(M,g)$ be a smooth, closed Riemannian manifold of dimension $m \geq 1$. Let $\lambda_2$ be the second lowest eigenvalue of $-\Delta_g$.  Let $A$ be a smooth, symmetric $(0,2)$-tensor field on $M$. Let $P_g$ be the differential operator on $M$ that is given on functions by $P_g u = \Delta^2_g  u - \text{div}_g(A(\nabla u)^\sharp)$ . Then, $P_g$ has a sign-changing  eigenfunction that is associated with the lowest eigenvalue of $P_g$, if $A$ has the property that $\max_{x \in M} \gamma_m(x) \leq -\lambda_2$ where $\gamma_m(x)$ is 
$$
\max_{w \in T_x(M) \setminus \{0\}}  \frac{A(x)(w,w)}{g(x)(w,w)}.
$$
($A(x)$ is  $A$ at $x \in M$, and $g(x)$ is $g$ at $x \in M$.) This, in turn, allows us to conclude that if $A$ has the property that $\max_{x \in M} \gamma_m(x) \leq -\lambda_2$, then $P_g$ is a fourth order Type 2 FSA elliptic differential operator.
\end{thm}

\begin{proof} First note that we have that
\begin{equation}
\int_M A((\nabla z)^\sharp,(\nabla z) ^\sharp) \, \text{dvol}_g \leq \int_M - \lambda_2 g((\nabla z)^\sharp,(\nabla z)^\sharp) \, \text{dvol}_g 
\end{equation}
for all functions $z \in H^2(M)$. Now, let $w$ be a non-constant solution of $-\Delta w = \lambda_2 w$. Then we have
\begin{equation}
\begin{split}
\int_M wP_gw \, \text{dvol}_g  = &  \int_M ((\Delta_g w)^2  + A((\nabla w)^\sharp ,(\nabla w)^\sharp)) \, \text{dvol}_g \\ \leq & \int_M ((\Delta_g w)^2 - \lambda_2 g((\nabla w)^\sharp, (\nabla w)^\sharp) \, \text{dvol}_g \\ = & \int_M ((\Delta_g w)^2 - \lambda_2 (-\Delta_g w)(w)) dvol_g \\ = & \int_M ({\lambda_2}^2 w^2 - {\lambda_2}^2 w^2) \, \text{dvol}_g \\ = & 0.
\end{split}
\end{equation}
It follows that either (a) there is a smooth function $z$ that has $\int_M zP_gz \, \text{dvol}_g < 0$, in which case Lemma 2.3 gives us the theorem, or (b) the lowest eigenvalue of $P_g$ is zero, and hence $w$, a sign-changing function, is an eigenfunction of $P_g$ associated with the lowest eigenvalue of $P_g$. The latter scenario  provides us with the theorem as well.
\end{proof}

\bibliographystyle{amsplain}

\begin{thebibliography}{9}



\bibitem{HR} Hebey, E., Robert, F. Coercivity and Struwe’s compactness for Paneitz fourth order operators with constant coefficients, Calculus of Variations and 

 Differential equations, 13, (2001), 491-517. 

\bibitem{Robert} Robert, F. Fourth order equations with critical growth in Riemannian geometry. Notes from lectures given at Madison and Berlin. https://iecl.univ-lorraine.fr/membre-iecl/robert-frederic-2/

\bibitem{Sweers} Sweers, G. On sign preservation for clotheslines, curtain rods, elastic membranes and thin plates, Jahresber. Dtsch. Math.-Ver. 118 (2016), 275–320. https://doi.org/10.1365/s13291-016-0147-0

\end{thebibliography}

\end{document}